\input amstex
\input amsppt.sty

\hsize13cm \vsize19cm \magnification1200

\TagsOnRight

\def\FKl{1}
\def\Ke{2}
\def\Kee{3}
\def\Ma1980{4}
\def\O1997{5}
\def\RB{6}
\def\Sa1996{7}
\def\StWhAB{8}
\def\StemAE{9}
\def\Vo{10}
\def\ZeilAG{11}

\def\AAa{1}
\def\AA{3}
\def\AB{4}
\def\AC{5}
\def\AD{6}
\def\AE{7}
\def\AF{8}
\def\AG{9}
\def\AH{10}
\def\AI{11}
\def\AJ{12}
\def\AK{13}
\def\AL{14}
\def\AM{2}

\def\TA{2}
\def\TB{3}
\def\TC{7}
\def\TD{1}

\def\CA{5}
\def\CB{6}

\def\({\left(}
\def\){\right)}
\def\be{\beta}
\def\la{\lambda}
\def\si{\sigma}
\def\ep{\varepsilon}
\def\et{\eta}

\def\S{\frak S}
\def\SS{{\frak S}}
\def\SP{\operatorname{SP}}
\def\Pf{\operatorname{Pf}}
\def\sgn{\operatorname{sgn}}
\def\tr{\operatorname{tr}}

\topmatter
\title Trace identities from identities for determinants
\endtitle
\author S.~Humphries and C.~Krattenthaler$^\dagger$
\endauthor
\affil Mathematics Department,
Brigham Young University\\
Provo, UT 84602, USA\\
e-mail: steve\@math.byu.edu\\
{WWW: \tt http://www.math.byu.edu/\~{}steve}\\\vskip6pt Institut
Girard
Desargues, Universit\'e Claude Bernard Lyon-I,\\
21, avenue Claude Bernard, F-69622 Villeurbanne Cedex, France.\\
e-mail: kratt\@euler.univ-lyon1.fr\\
WWW: \tt http://igd.univ-lyon1.fr/\~{}kratt
\endaffil
\address Mathematics Department,
Brigham Young University, Provo, UT 84602, USA. E-mail: {\tt
steve\@math.byu.edu}. {WWW: \tt
http://www.math.byu.edu/\~{}steve/}.
\endaddress
\address
Institut Girard Desargues, Universit\'e Claude Bernard Lyon-I, 21,
avenue Claude Bernard, F-69622 Villeurbanne Cedex, France. E-mail:
{\tt kratt\@euler.univ-lyon1.fr}.\linebreak WWW: \tt
http://igd.univ-lyon1.fr/\~{}kratt.
\endaddress
\thanks{$^\dagger$Research partially supported by EC's IHRP Programme,
grant HPRN-CT-2001-00272, ``Algebraic Combinatorics in
Europe"}\endthanks \subjclassyear{2000} \subjclass Primary 15A15;
Secondary 13B25 20G20
\endsubjclass
\keywords
\endkeywords
\abstract We present new identities for determinants of matrices
$(A_{i,j})$ with entries $A_{i,j}$ equal to $a_{i,j}$ or
$a_{i,0}a_{0,j}-a_{i,j}$, where the $a_{i,j}$'s are
indeterminates. We show that these identities are behind trace
identities for $SL(2,\Bbb C)$ matrices found earlier by Magnus in
his study of trace algebras.
\endabstract
\endtopmatter
\document

\leftheadtext{S. Humphries and C. Krattenthaler}


\subhead 1. Introduction\endsubhead In this paper we derive an
infinite family of new trace identities for $2 \times 2$ matrices
by using an infinite family of new determinantal identities. These
trace identities generalise a certain trace identity due to
Magnus.

Trace identities for $2 \times 2$ matrices have been studied for
over $100$ years, one of the original motivations being the
investigation of Teichm\" uller space via representations of
surface groups as (certain equivalence classes of) subgroups of
$SL(2,\Bbb C)$. This approach originated with Fricke and Klein
[\FKl] and there have been many subsequent attempts at ways of
giving real analytic trace coordinates for Teichm\" uller space
(see for example \cite{\Ke, \Kee, \O1997, \Sa1996} and references
therein).

Actions of groups on trace algebras have been investigated by Vogt
\cite{\Vo} and more recently by Magnus \cite{\Ma1980}. Vogt was
interested in studying invariants of differential equations, while
Magnus was concerned with automorphisms and outer automorphisms of
free groups. Physicists have also taken an interest in trace
relations \cite{\RB}.
  One thus sees the variety of applications that
these ideas have.

In \cite{\Ma1980}, Magnus's penultimate paper,  he investigated
the action of the automorphism group $Aut(F_n)$ of a free group
$F_n$ of rank $n$ on the traces of generic $2 \times 2$ matrices.
These generic traces generate
an algebra $Q_n$ as follows: let $m_1,\dots , m_n$ be `generic'
matrices in $SL(2,\Bbb Z)$ and for any sequence of distinct
elements  $i_1<i_2<\dots <i_k$, $k \le n$, of $\{1,2,\dots,n\}$ we
let
$$\tau_{i_1i_2\dots i_k}
(m_{1},m_{2},\dots, m_{n})=\tr (m_{i_1}m_{i_2}\cdots m_{i_k}),$$
where tr denotes the trace function. Then there are certain
relations among traces of $2 \times 2$ matrices that show that the
$\tau_i, \tau_{jk},\dots$ generate the
algebra of all the traces of elements of the group $\langle
m_1,\dots,m_n\rangle$: Take a polynomial algebra $Q_n'$ generated
by independent indeterminates $\tau_i', \tau_{jk}',\dots$; let $I
\subset Q_n'$ denote the ideal of all elements
$r(\tau_i',\tau_{jk}',\dots) \in Q_n'$ such that
$$r(\tau_i(m_1,\dots,m_n),\tau_{jk}(m_1,\dots,m_n),\dots) = 0$$ for
all choices of $m_i$. Then $Q_n=Q_n'/I$.

In \cite{\Ma1980} Magnus was concerned with obtaining
representations of $Out(F_n) = \mathbreak Aut(F_n)/Inn(F_n)$ by
looking at an induced action of $Out(F_n)$ on the trace
algebra $Q_n$. In order to do this he first investigates $Q_n$. In
doing this Magnus constructs various remarkable formulae satisfied
by the generators of the trace algebra. These are expressed as
equations in the determinants of certain matrices whose entries
are traces of elements of $SL(2,\Bbb C)$. He calls these the {\it
general identities} \cite{\Ma1980; p. 94} and uses them to derive
relations in the trace algebra that are needed for the proof that
$Q_n$ is finitely generated and other relevant properties.

These general identities are described as follows:

\proclaim {Magnus's Main Lemma} For $m_i, M_j \in SL(2,\Bbb C)$,
$1\le i, j \le 4$, we have the following trace relations:
$$\align &\det (\tr m_iM_j)+\det (\tr
m_iM_j^{-1})=0;\\&\det (\tr m_im_j)\det (\tr M_iM_j)=[\det (\tr
m_iM_j)]^2.\endalign$$\endproclaim

The goal of this paper is to search for the intrinsic background
of trace identities of the above kind. As a first observation, if
$m_i,M_j \in SL(2,\Bbb C)$, $0\le i,j$, with $m_0=M_0=I_2,$ the
$2\times 2$ identity matrix, and if we write $a_{i,j}=\tr m_iM_j$,
so that, in particular, $a_{i,0}=\tr m_i$ and $a_{0,j}=\tr M_j$
for $i,j\ge1$, then a simple calculation shows that
$$\tr m_iM_j^{-1}=a_{i,0}a_{0,j}-a_{i,j}.$$
Thus, we are led to search for identities involving determinants
of matrices $(A_{i,j})$, in which the entries $A_{i,j}$ may be of
the form $a_{i,j}$ or $a_{i,0}a_{0,j}-a_{i,j}$. We discovered that
on this abstract level there are in fact, somewhat surprisingly,
several of these. We summarise our findings in Theorems~\TD, \TB\
and \TC\ in the subsequent sections.

As a corollary to Theorem~\TD\ and to the trace theorem
Theorem~\TA, we obtain the following generalisation of the first
of Magnus's formulae.

\proclaim{Magnus's Main Lemma --- Generalised} For $n \ge 1$ let
$$m_1,\dots , m_n, M_1,\dots ,M_n \in SL(2,\Bbb C),$$ and put
$m_0=M_0=I_2.$ Define the $(n+1) \times (n+1)$ matrix
$A=(A_{i,j})_{0\le i,j\le n}$ by
$$A_{i,j}=\cases \tr(m_iM_j^{-1}) & \text {if }
i+j\text { is even;}\\
\tr(m_iM_j) &\text {otherwise},\endcases
$$
and define $n \times n$ matrices $B=(B_{i,j})_{1\le i,j\le n}$,
$C=(C_{i,j})_{1\le i,j\le n}$ by $B_{i,j}=-\tr m_iM_j, $ and $
C_{i,j}=\tr m_iM_j^{-1}$.
Then $$\det A = (-1)^n\det B + \det C.\tag\AAa$$ Further, if $n
\ge 4$, then $\det A=0$, while if $ n > 4$, then $\det B=\det
C=0$.\endproclaim

Theorem~\TD\ implies (\AAa), while the assertion in the last line
follows from Theorem~\TA.

Magnus's first identity is the above result in the situation
$n=4$, the first case where $\det A=0$, $\det B\ne0$ and $\det C
\ne 0.$.

In the next section, we state and prove Theorems~\TD\ and \TA. In
Section~3, we state and prove a general determinant formula, which
implies new trace identities for traces of the form $\tr m_im_j$
and $\tr m_im_j^{-1}$, where $m_1,m_2,\dots,m_n$ are given
matrices in $SL(2,\Bbb C)$. (That is, they address the case where
the second matrix family $M_1,M_2,\dots,M_n$ in our generalisation
of Magnus's Main Lemma is identical with the first one.)

Except for Theorem~\TA, which is not an abstract determinant
identity but a determinant identity specific for traces, we prove
our determinant identities by a combinatorial approach, as
proposed in \cite{\ZeilAG} (see also \cite{\StWhAB, Ch.~4}), that
is, we combinatorially expand both sides of our identities, and
then we bijectively identify the terms on the two sides, possibly
helped by an involution which cancels several terms on one side.

\example{Acknowledgement} Many thanks are due to Wayne Barrett for
connecting the authors.\endexample

\subhead 2. The determinant identities which imply the
generalisation of Magnus's formula\endsubhead In this section, we
prove a general determinant identity which implies our
generalisation (\AAa) of Magnus's formula (see Theorem~\TD\
below), and a general assertion about the vanishing of
determinants formed out of traces (see Theorem~\TA\ below) that
implies the last assertion in our generalisation of Magnus's Main
Lemma, but, in addition, produces many more trace identities.

\proclaim{Theorem \TD} Let $(a_{i,j})_{0\le i,j\le n}$ be a doubly
indexed sequence with the property that $a_{0,0}=2$. We let $A$ be
the $(n+1)\times (n+1)$ matrix $A=(A_{i,j})_{1\le i,j\le n}$,
where
$$A_{i,j}=\cases
\la a_{0,j}&\text {if }0=i\ne j,\\
a_{i,j}&\text {if }j=0\text{ or if }i+j\text { is even},\\
\la a_{i,0}a_{0,j}-a_{i,j}&\text {otherwise.}\\
\endcases$$
Furthermore, we define two $n\times n$ matrices $B=(B_{i,j})_{1\le
i,j\le n}$ and $C=(C_{i,j})_{1\le i,j\le n}$ by
$$B_{i,j}=
\la a_{i,0}a_{0,j}-a_{i,j},$$ and
$$C_{i,j}=
a_{i,j}.$$ Then
$$\det A-(-1)^n\det B-\det C=0.\tag\AM$$
\endproclaim

\remark{Remark} For better clarity, we remark that, by our
definitions, the first row of $A$ is
$$(2,\la a_{0,1},\la a_{0,2},\dots,\la a_{0,n}),$$
while the first column is
$$(2, a_{1,0}, a_{2,0},\dots, a_{n,0})^t.$$
For example, for $n=4$, the matrix $A$ is equal to
$$\pmatrix
        2 & \la a_{0,1} & \la a_{0,2} &
           \la a_{0,3} & \la a_{0,4}\\
         a_{1,0} & a_{1,1} & \la a_{1,0} a_{0,2} - a_{1,2} & a_{1,3} &
           \la a_{1,0} a_{0,4} - a_{1,4}\\
         a_{2,0} & \la a_{2,0} a_{0,1} - a_{2,1} & a_{2,2} &
           \la a_{2,0} a_{0,3} - a_{2,3} & a_{2,4}\\
         a_{3,0} & a_{3,1} & \la a_{3,0} a_{0,2} - a_{3,2} & a_{3,3} &
           \la a_{3,0} a_{0,4} - a_{3,4}\\
         a_{4,0} & \la a_{4,0} a_{0,1} - a_{4,1} & a_{4,2} &
           \la a_{4,0} a_{0,3} - a_{4,3} & a_{4,4}\\
        \endpmatrix.
$$
In view of this remark, it should be clear that the case $\la=1$
of this theorem implies (\AAa).
\endremark

\demo{Proof} The first observation is that, in $\det A$, the
coefficient of $\la^m$ is zero for $m\ge2$, and the same is true
for $\det B$ and $\det C$. Clearly, this is trivial for $\det C$,
which contains no $\la$ at all. To see the claim for $\det B$, we
let $B'$ be the matrix $(\la a_{i,0}a_{0,j})_{1\le i,j\le n}$ and
$B''$ the matrix $(-a_{i,j})_{1\le i,j\le n}$. Then we expand
$\det B=\det(B'+B'')$ in the following form,
$$\det B=
\sum _{J\subseteq \{1,2,\dots,n\}} ^{}\det B^J,$$ where $B^J$
denotes the matrix where the columns indexed by elements from $J$
are those from $B'$, while the remaining columns are those from
$B''$. Since any two columns from $B'$ are dependent, we have
$\det B^J=0$ whenever $\vert J\vert\ge2$. Thus the largest
exponent of $\la$ in the expansion of $\det B$ is 1.

For $\det A$ we proceed in the same way. We let $A'$ be the matrix
which contains the ``$\la$-terms" from $A$, and we let $A''$ be
the ``rest". To be precise,
$$A'_{i,j}=\cases
\la a_{0,j}&\text {if }0=i\ne j,\\
0&\text {if }j=0\text{ or if }i+j\text { is even},\\
\la a_{i,0}a_{0,j}&\text {otherwise.}\\
\endcases,$$
while
$$A''_{i,j}=\cases
0&\text {if }0=i\ne j,\\
a_{i,j}&\text {if }j=0\text{ or if }i+j\text { is even},\\
-a_{i,j}&\text {otherwise.}\\
\endcases$$
Again, we have chosen $A'$ and $A''$ so that $A=A'+A''$. Then we
do the same expansion as before,
$$\det A=
\sum _{J\subseteq \{1,2,\dots,n\}} ^{}\det A^J,$$ with the
analogous meaning of $A^J$. Again, we have $\det A^J=0$ if $\vert
J\vert\ge2$, this time because the $0$-th column of $A'$ is 0, and
because any two columns of $A'$ indexed by $j_1,j_2\ge1$ which
have the same parity are dependent, while if $j_1$ and $j_2$ have
different parity, $1/\la a_{0,j_1}$ times the $j_1$-st column plus
$1/\la a_{0,j_2}$ times the $j_2$-nd column gives the first column
of $A''$.

It remains to verify that the coefficients of $\la^0$ and of
$\la^1$ in (\AM) vanish.

Let us begin with the coefficient of $\la^0$. If we set $\la=0$ in
(\AM), then $\det A$ can be reduced to
$$2\det ((-1)^{i+j}a_{i,j})_{1\le i,j\le n}=
2\det (a_{i,j})_{1\le i,j\le n}=2\det C,$$ while
$$\det B=\det (-a_{i,j})_{1\le i,j\le n}=
(-1)^n\det C.$$ Thus, the coefficient of $\la^0$ in (\AM) is
indeed zero.

For the coefficient of $\la^1$ we only have to look at $\det A$
and $\det B$. We shall derive combinatorial expressions for these
two determinants. In order to do so, let us write $\S_n$ for the
symmetric group on $\{1,2,\dots, n\}$. From the definition of the
determinant, we have
$$\det B=\sum _{\si\in\S_n} ^{}\sgn \si\prod _{i=1} ^{n}B_{i,\si(i)}.
$$
Extracting the coefficient of $\la^1$, we see that the coefficient
of $\la^1$ in $\det B$ is
$$\align
{\sum _{\si\in\S_n} ^{}}\sgn \si \sum _{k=1} ^{n}
a_{k,0}a_{0,\si(k)}&
\underset i\ne k\to{\prod _{i=1} ^{n}}(-a_{i,\si(i)})\\
&= \sum _{k=1} ^{n} \sum _{l=1} ^{n} (-1)^{k+l} a_{k,0}a_{0,l}
{\sum _{\si\in\S_n^{(k,l)}} ^{}} \sgn \si
\underset i\ne k\to{\prod _{i=1} ^{n}}(-a_{i,\si(i)})\\
&=(-1)^{n-1} \sum _{k=1} ^{n} \sum _{l=1} ^{n} (-1)^{k+l}
a_{k,0}a_{0,l} {\sum _{\si\in\S_n^{(k,l)}} ^{}} \sgn \si \underset
i\ne k\to{\prod _{i=1} ^{n}}a_{i,\si(i)},
\endalign$$
where $\S_n^{(k,l)}$ is the set of bijections from
$\{1,2,\dots,n\} \backslash \{k\}$ to $\{1,2,\dots,n\} \backslash
\{l\}$, where $\sgn\si$ has the obvious meaning when identifying
$\S_n^{(k,l)}$ with $\S_{n-1}$.

Aplying a similar procedure to $\det A$, we obtain that the
coefficient of $\la^1$ in $\det A$ is
$$\align
2 {\sum _{\si\in\S_n} ^{}} &\underset k+\si(k)\text{ odd}\to{\sum
_{k=1} ^{n}} (\sgn \si)\cdot a_{k,0}a_{0,\si(k)}
\underset i\ne k\to{\prod _{i=1} ^{n}}(-1)^{i+\si(i)}a_{i,\si(i)}\\
&\kern2cm + \sum _{k=1} ^{n} \sum _{l=1} ^{n}(-1)^{k+l-1}
a_{k,0}a_{0,l} {\sum _{\si\in\S_n^{(k,l)}} ^{}}\sgn \si
\underset i\ne k\to{\prod _{i=1} ^{n}}(-1)^{i+\si(i)}a_{i,\si(i)}\\
&=2\underset k+l\text{ odd}\to{\sum _{k,l=1} ^{n}} a_{k,0}a_{0,l}
{\sum _{\si\in\S_n^{(k,l)}} ^{}} \sgn \si \underset i\ne
k\to{\prod _{i=1} ^{n}}a_{i,\si(i)} - \sum _{k=1} ^{n} \sum _{l=1}
^{n}a_{k,0}a_{0,l} {\sum _{\si\in\S_n^{(k,l)}} ^{}}\sgn \si
\underset i\ne k\to{\prod _{i=1} ^{n}}a_{i,\si(i)}\\
&= \sum _{k=1} ^{n} \sum _{l=1} ^{n}(-1)^{k+l-1}a_{k,0}a_{0,l}
{\sum _{\si\in\S_n^{(k,l)}} ^{}}\sgn \si \underset i\ne k\to{\prod
_{i=1} ^{n}}a_{i,\si(i)}.
\endalign
$$
Thus, indeed, the coefficient of $\la^1$ in (\AM) is zero.

This completes the proof of the theorem. \qed
\enddemo


The last assertion of our generalisation of Magnus's Main Lemma
follows from the following more general result which gives more
trace relations.

\proclaim{Theorem \TA} Let $ m_1,m_2,\dots , m_n,M_1,M_2,\dots,
M_n \in SL(2,\Bbb C)$ and let $
\varepsilon_1,\varepsilon_2,\dots,\varepsilon_n \in \{\pm 1\}$.
Define the $n \times n$ matrix $D=(D_{i,j})_{1\le i,j\le n}$ by
$D_{i,j}=\tr (m_iM_j^{\varepsilon_i}).$ If $n \ge 5$, then $\det D
= 0.$\endproclaim

\demo{Proof} The proof will follow by exhibiting a $0$-eigenvector
for $D$.

Fix $j \le n$ and let $M_j=\pmatrix
M_{j1}&M_{j2}\\M_{j3}&M_{j4}\endpmatrix$, so that
$M_j^{-1}=\pmatrix\hphantom{-}M_{j4}&-M_{j2}\\
-M_{j3}&\hphantom{-}M_{j1}\endpmatrix$. We also let $m_i=\pmatrix
m_{i1}&m_{i2}\\m_{i3}&m_{i4}\endpmatrix$ for all $i \le n$.  We
will find non-trivial functions $v_1,\dots,v_n$ of the variables
$m_{i,j}$ such that $v=(v_1,\dots, v_n)$ is a (left)
$0$-eigenvector for $D$, so that $vD=0$. This will be the case if
for all $j \le n$ we have $$\sum_{i=1}^n v_i\tr
(m_iM_j^{\varepsilon_i})=0.$$ But each of the above equations is
linear  in the variables $M_{j1}, M_{j2}, M_{j3}, M_{j4}$ and the
equations that we so obtain are independent of the column index
$j$. We thus have $n
> 4$ linear equations in $4$ unknowns. There is thus, generically,
a non-trivial solution. Thus, $\det D=0$ except on a set of
measure zero; but since $\det$ is continuous it follows that we
always have $\det D=0$.\qed
\enddemo

\subhead 3. More determinant identities\endsubhead Here, we
present an identity for determinants of matrices, in which the
entries are $a_{i,j}$ or $a_ia_j-a_{i,j}$. Interpreting $a_i$ as
the trace of a matrix $m_i\in SL(2,\Bbb C)$ and $a_{i,j}$ as $\tr
m_im_j$, these identities produce therefore trace identities for
the traces $\tr m_im_j$ and $\tr m_im_j^{-1}$. In fact, in the
theorem below, we allow two additional parameters, $\la$ and
$\be$. By specialising them in different ways, we obtain various
new determinant identities. The specialisation in Corollary~\CA\
contains the identity which is relevant to the trace case, whereas
Corollary~\CB\ contains a ``skew" variation. (The ``skew" refers
to the fact that the matrix $A$ there is skew-symmetric.) As an
aside, we prove in Theorem~\TC\ the curious fact, that, in the
``skew" case, the determinant of the matrix $C$ factors into two
big factors, one of which collects the ``even terms" of the
Pfaffian of $A$, the other collecting its odd terms.

We alert the reader that, when compared to Theorem~\TD, the two
theorems below follow a different index convention in that the
entries of the matrix $A$ are indexed by $i$ and $j$ from
$\{1,2,\dots,n\}$ (rather than $\{0,1,\dots,n\}$), and, similarly,
the entries of the matrices $B$ and $C$ are indexed by $i$ and $j$
from $\{2,3,\dots,n\}$ (rather than $\{1,2,\dots,n\}$). This
convention has advantages over the other in the formulation of
Theorem~\TC. A further change of convention is that the
``trace-like" entries are $a_{i,j}-\la a_{1,i}a_{1,j}$ (rather
than $\la a_{1,i}a_{1,j}-a_{i,j}$). This allows a more elegant
formulation of the following theorem, but, clearly, by multiplying
every other row and column of $B$ and $C$ by $-1$, we could pass
to the convention for the ``trace-like" entries which is followed
in Theorem~\TD.

\proclaim{Theorem \TB} Let $(a_{i,j})_{1\le i,j\le n}$ be a doubly
indexed sequence with the property that  $a_{i,1}=\be a_{1,i}$ for
$i>1$. We let $A$ be the $n\times n$ matrix $A=(A_{i,j})_{1\le
i,j\le n}$, where
$$A_{i,j}=\cases
\la a_{1,j}&\text {if }1=i\ne j,\\
a_{i,j}&\text{otherwise,}
\endcases$$
Furthermore, we define two $(n-1)\times (n-1)$ matrices
$B=(B_{i,j})_{2\le i,j\le n}$ and $C=(C_{i,j})_{2\le i,j\le n}$ by
$$B_{i,j}=\cases
a_{i,j}-\la a_{1,i}a_{1,j}&\text {if }i+j\text { is even,}\\
a_{i,j}&\text {if }i+j\text { is odd},\endcases$$ and
$$C_{i,j}=\cases
a_{i,j}&\text {if }i+j\text { is even},\\
a_{i,j}-\la a_{1,i}a_{1,j}&\text {if }i+j\text { is odd.}\\
\endcases$$
Then
$$\det A-\be(\det B+\det C)=(a_{1,1}-2\be)
\det(a_{i,j})_{2\le i,j\le n}.\tag\AA$$
\endproclaim

\remark{Remark} For the benefit of the reader, we display the
matrices $A$, $B$, and $C$ for $n=5$:
$$\gather
A= \pmatrix\format\r&\quad \r&\quad \r&\quad \r&\quad \r\\
         a_{1,1} & \la a_{1,2} & \la a_{1,3} &
            \la a_{1,4} & \la a_{1,5}\\
         \be a_{1,2} & a_{2,2} & a_{2,3} & a_{2,4} & a_{2,5}\\
         \be a_{1,3} & a_{3,2} & a_{3,3} & a_{3,4} & a_{3,5}\\
         \be a_{1,4} & a_{4,2} & a_{4,3} & a_{4,4} & a_{4,5}\\
         \be a_{1,5} & a_{5,2} & a_{5,3} & a_{5,4} & a_{5,5}\\
         \endpmatrix,\\
B= \pmatrix
          a_{2,2} - \la a_{1,2}^2   & a_{2,3} &
            - \la a_{1,2} a_{1,4}   + a_{2,4} &
            a_{2,5}\\
         a_{3,2} &   a_{3,3}- \la a_{1,3}^2   &
            a_{3,4} &a_{3,5} - \la a_{1,3} a_{1,5}
            \\
          a_{4,2} - \la a_{1,2} a_{1,4}   & a_{4,3} &
            a_{4,4}  - \la a_{1,4}^2   & a_{4,5}\\
         a_{5,2} &  a_{5,3}- \la a_{1,3} a_{1,5}    &
            a_{5,4} &  a_{5,5}- \la a_{1,5}^2   \\
         \endpmatrix,\\
C= \pmatrix
         a_{2,2} &  a_{2,3}- \la a_{1,2} a_{1,3}    &
            a_{2,4} & a_{2,5}- \la a_{1,2} a_{1,5}
            \\
           a_{3,2}- \la a_{1,2} a_{1,3}   & a_{3,3} &
             a_{3,4}- \la a_{1,3} a_{1,4}    &
            a_{3,5}\\
         a_{4,2} &  a_{4,3}- \la a_{1,3} a_{1,4}    &
            a_{4,4} & a_{4,5}- \la a_{1,4} a_{1,5}
            \\
           a_{5,2}- \la a_{1,2} a_{1,5}   & a_{5,3} &
             a_{5,4} - \la a_{1,4} a_{1,5}   &
            a_{5,5}\\
         \endpmatrix.
\endgather$$
\endremark

\demo{Proof} As a first step, we derive combinatorial expressions
for $\det A$, $\det B$ and $\det C$.

Let us introduce some notation. As earlier, we write $\S_n$ for
the symmetric group on $\{1,2,\dots, n\}$. Given a permutation
$\si\in\S_n$, we let $c(\si)$ be the number of cycles of $\si$.
Furthermore we define the function $f_1$ by
$$f_1(\si)=\cases 1&\text {if }\si(1)=1,\\
0&\text {if }\si(1)\ne1.\endcases$$

By the definition of the determinant, we have
$$\det A=\sum _{\si\in\S_n} ^{}\sgn \si\prod _{i=1} ^{n}A_{i,\si(i)}=
\sum _{\si\in\S_n} ^{}(-1)^{n-c(\si)} \la^{1-f_1(\si)} \prod
_{i=1} ^{n}a_{i,\si(i)}. \tag\AB$$

\medskip
In order to describe combinatorial expressions for $\det B$ and
$\det C$, we need some further notations and definitions. We write
$\SS_{n-1}$ for the symmetric group on\linebreak $\{2,3,\dots,
n\}$. A {\it signed\/} permutation on $\{2,3,\dots,n\}$ is a pair
$(\pi,\ep)$, where $\pi\in\SS_{n-1}$ and $\ep\in \{-1,1\}^{n-1}$.
For the sake of convenience, we label the components of $\ep$ from
2 through $n$, that is, $\ep=(\ep_2,\ep_3,\dots,\ep_n)$. We define
the {\it weight\/} $w(\pi,\ep)$ of a signed permutation
$(\pi,\ep)$ by
$$w(\pi,\ep)=\prod _{i=2} ^{n}W_{i,\pi(i)},$$
where
$$W_{i,\pi(i)}=\cases a_{i,\pi(i)}&\text {if }\ep_i=1,\\
-\la a_{1,i}a_{1,\pi(i)}&\text {if }\ep_i=-1.\endcases$$ We need
two particular subsets of all signed permutations: let
$\SP_{n-1}^{(1)}$ denote the set of all signed permutations on
$\{2,3,\dots,n\}$ with $\ep_i=1$ whenever $i+\pi(i)$ is odd, $2\le
i\le n$, and let $\SP_{n-1}^{(2)}$ denote the set of all signed
permutations on $\{2,3,\dots,n\}$ with $\ep_i=1$ whenever
$i+\pi(i)$ is even, $2\le i\le n$.

Now, with all this notation, for the determinant of $B$ we have
$$\det B=\sum _{\si\in\SS_{n-1}} ^{}\sgn \si\prod _{i=2} ^{n}B_{i,\si(i)}=
\sum _{(\pi,\ep)\in \SP^{(1)}_{n-1}} ^{}(\sgn \pi)\cdot
w(\pi,\ep), \tag\AC$$ while for the determinant of $C$ we have
$$\det C=\sum _{\si\in\SS_{n-1}} ^{}\sgn \si\prod _{i=2} ^{n}B_{i,\si(i)}=
\sum _{(\pi,\ep)\in \SP^{(2)}_{n-1}} ^{}(\sgn \pi)\cdot
w(\pi,\ep). \tag\AD$$

We can now begin the actual proof of (\AA).

We start by identifying the coefficients of $\la^0$: By
inspection, the coefficient of $\la^0$ in $\det A$ is
$a_{1,1}\det(a_{i,j})_{2\le i,j\le n}$. It is equally obvious that
the coefficient of $\la^0$ in $\det B$, as well as in $\det C$, is
equal to $\det(a_{i,j})_{2\le i,j\le n}$. Thus, the coefficients
of $\la^0$ on both sides of (\AA) agree.

Next we identify the coefficients of $\la^1$: clearly, the
coefficient of $\la^1$ on the right-hand side of (\AA) is zero. We
now use expressions (\AB)--(\AD) for $\det A$, $\det B$, and $\det
C$, respectively, to show that this is also the case on the
left-hand side. The coefficient of $\la^1$ in (\AB) is
$$\align
\underset \si(1)\ne1\to{\sum _{\si\in\S_n} ^{}} (-1)^{n-c(\si)}
\prod _{i=1} ^{n}a_{i,\si(i)} &= \underset \si(1)\ne1\to{\sum
_{\si\in\S_n} ^{}} (-1)^{n-c(\si)}
a_{\si^{-1}(1),1}a_{1,\si(1)}\underset i\ne \si^{-1}(1)\to {\prod
_{i=2}
^{n}}a_{i,\si(i)}\\
&= \be\underset \si(1)\ne1\to{\sum _{\si\in\S_n} ^{}}
(-1)^{n-c(\si)} a_{1,\si^{-1}(1)}a_{1,\si(1)}\underset i\ne
\si^{-1}(1)\to {\prod _{i=2} ^{n}}a_{i,\si(i)}. \tag\AE\endalign$$
Terms contributing to the coefficient of $\la^1$ in (\AC) and
(\AD) occur exactly for the signed permutations $(\pi,\ep)$ where
$\ep$ is a vector with exactly one component equal to $-1$. Let
$\ep^{(k)}$ denote the vector with all components equal to 1
except for the $k$-th, which is $-1$. The coefficient of $\la^1$
in (\AC) is then
$$
\sum _{\pi\in \SS_{n-1}} ^{} \underset k+\pi(k)\text {
even}\to{\sum _{k=2} ^{n}} (\sgn \pi)\cdot \frac
{w(\pi,\ep^{(k)})}\la= -\kern-4pt\sum _{\pi\in \SS_{n-1}} ^{}
\underset k+\pi(k)\text { even}\to{\sum _{k=2} ^{n}} (\sgn
\pi)\cdot a_{1,k}a_{1,\pi(k)} \underset i\ne k\to{\prod _{i=2}
^{n}}a_{i,\pi(i)},
$$
while the coefficient of $\la^1$ in (\AD) is
$$
\sum _{\pi\in \SS_{n-1}} ^{} \underset k+\pi(k)\text {
odd}\to{\sum _{k=2} ^{n}} (\sgn \pi)\cdot \frac
{w(\pi,\ep^{(k)})}\la= -\sum _{\pi\in \SS_{n-1}} ^{} \underset
k+\pi(k)\text { odd}\to{\sum _{k=2} ^{n}} (\sgn \pi)\cdot
a_{1,k}a_{1,\pi(k)} \underset i\ne k\to{\prod _{i=2}
^{n}}a_{i,\pi(i)}. \tag\AF$$ The coefficient of $\la^1$ in the sum
of (\AC) and (\AD) is therefore
$$
-\sum _{\pi\in \SS_{n-1}} ^{} {\sum _{k=2} ^{n}}
(-1)^{n-1-c(\pi)}\cdot a_{1,k}a_{1,\pi(k)} \underset i\ne
k\to{\prod _{i=2} ^{n}}a_{i,\pi(i)}. \tag\AG$$ To a fixed pair
$(\pi,k)$, $\pi\in\SS_{n-1}$ and $2\le k\le n$, we now associate
the permutation $\si\in\S_n$, given by $\si(k)=1$,
$\si(1)=\pi(k)$, and $\si(i)=\pi(i)$ for all $i\ne 1,k$. Thus we
see that (\AE) and (\AG) multiplied by $\be$ are identical, which
implies that the coefficient of $\la^1$ on the left-hand side of
(\AA) vanishes, as desired.

The remaining task is to show that all other terms in the sum of
(\AC) and (\AD) cancel. The reader should note that these ``other
terms" in (\AC) and (\AD) are indexed by signed permutations
$(\pi,\ep)$ in the union $\SP^{(1)}_{n-1}\cup \SP^{(2)}_{n-1}$,
where the vector $\ep$ has at {\it least two} components equal to
$-1$. We show that these terms cancel by defining a {\it
sign-reversing involution} $i$ on the set of signed permutations
in $\SP^{(1)}_{n-1}\cup \SP^{(2)}_{n-1}$, where the vector $\ep$
has at least two components equal to $-1$. The map $i$ will be
defined separately on three disjoint subsets of this set.

\smallskip
{\it Set 1}. Consider all signed permutations $(\pi,\ep)$ in
$\SP^{(1)}_{n-1}\cup \SP^{(2)}_{n-1}$ with the property that there
are at least two {\it even} indices $i_1,i_2$, $2\le i_1<i_2\le
n$, with $\ep_{i_1}=\ep_{i_2}=-1$. Let us call this property the
property~P1.

Given a signed permutation $(\pi,\ep)$ with property~P1, let $i_1$
and $i_2$ be even integers such that $\ep_{i_1}=\ep_{i_2}=-1$,
which are minimal with respect to this property. Then we define
$$i\big((\pi,\ep)\big):=(\pi\circ (i_1,i_2),\ep).\tag\AH$$
(The permutation $\pi\circ (i_1,i_2)$ is the composition of $\pi$
and the transposition exchanging $i_1$ and $i_2$.) Clearly,
$i\big((\pi,\ep)\big)$ also has property~P1 since the vector $\ep$
has not changed. Furthermore, if $(\pi,\ep)\in\SP^{(1)}_{n-1},$
then also $i\big((\pi,\ep)\big)\in\SP^{(1)}_{n-1}$, and similarly
for $\SP^{(2)}_{n-1}$. By definition, the weight of $(\pi,\ep)$ is
$$w(\pi,\ep)=(-\la a_{1,i_1}a_{1,\pi(i_1)})(-\la
a_{1,i_2}a_{1,\pi(i_2)}) \underset i\ne i_1,i_2\to{\prod _{i=2}
^{n}}W_{i,\pi(i)},$$ while the weight of $i\big((\pi,\ep)\big)$ is
$$w\big(\pi\circ (i_1,i_2),\ep\big)=(-\la a_{1,i_1}a_{1,\pi(i_2)})(-\la
a_{1,i_2}a_{1,\pi(i_1)}) \underset i\ne i_1,i_2\to{\prod _{i=2}
^{n}}W_{i,\pi(i)},$$ which is equal to $w(\pi,\ep)$. In summary,
we have established the relation
$$(\sgn \pi)\cdot w(\pi,\ep)=
-(\sgn (\pi\circ(i_1,i_2)))\cdot w\big(\pi\circ
(i_1,i_2),\ep\big).$$ Since, in addition, $i$ is an involution,
the terms in (\AC) indexed by signed permutations with property~P1
cancel each other pairwise, and the same is true for the analogous
terms in (\AD).

\smallskip
{\it Set 2}. Now we consider all signed permutations in
$\SP^{(1)}_{n-1}\cup \SP^{(2)}_{n-1}$ which do not have
property~P1, but have the property that there are at least two
{\it odd} indices $i_1,i_2$, $2\le i_1<i_2\le n$ with
$\ep_{i_1}=\ep_{i_2}=-1$. Let us call this property the
property~P2. Given a signed permutation $(\pi,\ep)$ with
property~P2, we define the map $i$ by (\AH), as before. It is then
easy to see that everything else is as in the previous case. In
particular, the terms in (\AC) indexed by signed permutations with
property~P2 cancel each other pairwise, and the same is true for
the analogous terms in (\AD).

\smallskip
{\it Set 3}. Finally, we consider the signed permutations
$(\pi,\ep)$ in $\SP^{(1)}_{n-1}\cup \SP^{(2)}_{n-1}$ which have
neither property~P1 nor property~P2. Since $\ep$ must have at
least two components equal to $-1$, the only possibility is then
that there is an {\it even} $i_1$ and an {\it odd} $i_2$ such that
$\ep_{i_1}=\ep_{i_2}=-1$, and these are the only components of
$\ep$ which are equal to $-1$. We define the map $i$ again by
(\AH). This time, if $(\pi,\ep)\in\SP^{(1)}_{n-1}$, then
$i\big((\pi,\ep)\big)\in\SP^{(2)}_{n-1}$, and vice versa. However,
all the other conclusions of the first case remain valid, and,
thus, again the terms in the sum of the right-hand sides of (\AC)
and (\AD) indexed by signed permutations in this subset cancel
each other pairwise.

\medskip
This completes the proof of (\AA). \qed
\enddemo

\proclaim{Corollary \CA}Let $(a_{i,j})_{1\le i,j\le n}$ be the
doubly indexed sequence with the property that $a_{i,j}=a_{j,i}$
and $a_{i,i}=2$ for all $i$ and $j$. Then, if the matrices $A$,
$B$ and $C$ are defined as in Theorem~\TB\ {\rm(}with
$\be=1${\rm)}, we have
$$\det A=\det B+\det C.$$
\endproclaim
\demo{Proof}We set $\be=1$, $a_{i,j}=a_{j,i}$ and $a_{i,i}=2$ for
all $i$ and $j$ in Theorem~\TB. Then $a_{1,1}-2\be=0$, and, hence,
the assertion is equivalent to Equation~(\AA) with these
specializations.\qed
\enddemo

\proclaim{Corollary \CB}Let $(a_{i,j})_{1\le i,j\le n}$ be the
doubly indexed sequence with the property that $a_{i,j}=-a_{j,i}$
for all $i$ and $j$. Then, if the matrices $A$, $B$ and $C$ are
defined as in Theorem~\TB\ {\rm(}with $\be=-1${\rm)}, for even $n$
we have
$$\det A+\det B+\det C=0.$$
\endproclaim

\demo{Proof}We set $\be=-1$ and $a_{i,j}=-a_{j,i}$ for all $i$ and
$j$ in Theorem~\TB. Then the determinant on the right-hand side of
(\AA) is the determinant of a skew-symmetric matrix of odd size
and, hence, zero. The assertion is therefore equivalent to
Equation~(\AA) with these specializations.\qed
\enddemo

As it turns out, in the skew-symmetric case, the determinant of
the matrix $C$ factors into two big factors. These two factors can
be described explicitly. They are the ``even" and the ``odd" parts
of the Pfaffian of $A$. To be more precise, we put ourselves in
the special case of Theorem~\TB where $n$ is even, $\la=1$, and
$a_{i,j}=-a_{j,i}$ for all $i$ and $j$ (in particular, $\be=-1$).
By definition, the Pfaffian of $A$, denoted $\Pf(A)$, is the
square root of $\det A$, where the sign is chosen so that the term
$a_{1,n}a_{2,n-1}\cdots a_{n/2,n/2+1}$ occurs with coefficient
$+1$. In combinatorial terms, the Pfaffian of $A$ is (cf\.
\cite{\StemAE, Sec.~2})
$$\Pf(A)=
\sum _{m} ^{}\sgn m\cdot \prod _{(i,j)\in m} ^{}a_{i,j},\tag\AI$$
where the sum is over all perfect matchings $m$ on
$\{1,2,\dots,n\}$. Let $\Pf_e(A)$ denote the sum of all the terms
appearing on the right-hand side of (\AI) which contain $a_{1,k}$
for an {\it even} $k$. Similarly, we denote by $\Pf_o(A)$  the sum
of all the terms appearing in $\Pf(A)$ which contain $a_{1,k}$ for
an {\it odd\/} $k$.

\proclaim{Theorem \TC} Let $(a_{i,j})_{1\le i,j\le n}$ be the
doubly indexed sequence with the property that $a_{i,j}=-a_{j,i}$
for all $i$ and $j$. Then, if the matrices $A$ and $C$ are defined
as in Theorem~\TB\ {\rm(}with $\be=-1${\rm)}, for even $n$ we have
$$\det C=-2\Pf_e(A)\Pf_o(A).$$
\endproclaim
\demo{Proof} We expand the determinant $\det(C)$ as in the proof
of Theorem~\TB, see (\AD). Clearly, in this expansion, all the
contributions of signed permutations in $\SP^{(2)}_{n-1}$ for
which all the $\ep_i$'s are 1 cancel each other, because the sum
of these contributions is simply $\det(a_{i,j})_{2\le i,j\le n}$,
the determinant of a skew-symmetric matrix of odd dimension. The
arguments given in the proof of Theorem~\TB\ for Sets~1 and 2 show
in addition that the contributions of signed permutations in
$\SP^{(2)}_{n-1}$ for which there exist $i_1$ and $i_2$, both even
or both odd, such that $\ep_{i_1}=\ep_{i_2}=-1$ cancel each other.
However, this happens also for signed permutations in Set~3, i.e.,
for signed permutations in $\SP^{(2)}_{n-1}$ for which
$\ep_{i_1}=\ep_{i_2}=-1$ for an even $i_1$ and an odd $i_2$, and
for which all other $\ep_i$'s are 1. To see this, for fixed even
$i_1$ and odd $i_2$, let us consider all the signed permutations
$(\pi,\ep)$ in  $\SP^{(2)}_{n-1}$ with $\ep_{i_1}=\ep_{i_2}=-1$,
$\ep_i=1$ otherwise, with a fixed value $\pi(i_1)$ and a fixed
value $\pi(i_2)$. Their contribution to (\AD) is
$$\sgn(i_1,i_2,\pi(i_1),\pi(i_2))\cdot a_{1,i_1}a_{1,\pi(i_1)}
a_{1,i_2}a_{1,\pi(i_2)}\det
A^{1,i_1,i_2}_{1,\pi(i_1),\pi(i_2)},\tag\AJ$$ where
$A^{1,i_1,i_2}_{1,\pi(i_1),\pi(i_2)}$ is the submatrix of $A$
which arises by deleting the rows numbered $1$, $i_1$ and $i_2$
and the columns numbered $1$, $\pi(i_1)$ and $\pi(i_2)$, and
where\linebreak $\sgn(i_1,i_2,\pi(i_1),\pi(i_2))$ is a certain
sign. If we interchange the roles of $i_1$ and $\pi(i_1)$, and of
$i_2$ and $\pi(i_2)$, and consider the analogous signed
permutations in $\SP^{(2)}_{n-1}$, then their contribution is
$$\sgn(i_1,i_2,\pi(i_1),\pi(i_2))\cdot a_{1,i_1}a_{1,\pi(i_1)}
a_{1,i_2}a_{1,\pi(i_2)}\det
A_{1,i_1,i_2}^{1,\pi(i_1),\pi(i_2)}.\tag\AK$$ However, we have
$$A^{1,i_1,i_2}_{1,\pi(i_1),\pi(i_2)}=
-\left(A_{1,i_1,i_2}^{1,\pi(i_1),\pi(i_2)}\right)^t,$$ which
implies
$$\det A^{1,i_1,i_2}_{1,\pi(i_1),\pi(i_2)}=
-\det A_{1,i_1,i_2}^{1,\pi(i_1),\pi(i_2)},$$ since these are
determinants of matrices of dimension $n-3$, which is odd. Thus,
the sum of the terms (\AJ) and (\AK) is zero.

In summary, the above arguments have shown that $\det C$ is equal
to the contributions in (\AD) by signed permutations in
$\SP^{(2)}_{n-1}$ with {\it exactly} one $\ep_i$ which is $-1$. To
be precise, they show that (compare with the expression (\AF))
$$\det C=-\sum _{\pi\in \SS_{n-1}} ^{}
\underset k+\pi(k)\text { odd}\to{\sum _{k=2} ^{n}} (\sgn
\pi)\cdot a_{1,k}a_{1,\pi(k)} \underset i\ne k\to{\prod _{i=2}
^{n}}a_{i,\pi(i)}.$$

In fact, there is still cancellation in the expression on the
right-hand side. If $\pi\in \SS_{n-1}$ should have a cycle of odd
length which does not contain $k$, then the permutation $\bar \pi$
arising from $\pi$ by reversing the orientation of the cycle has
the same sign, but the product ${\prod _{i=2,i\ne k}
^{n}}a_{i,\bar\pi(i)}$ has sign opposite to ${\prod _{i=2,i\ne k}
^{n}}a_{i,\pi(i)}$. Thus, the contributions corresponding to $\pi$
and $\bar\pi$ cancel each other. An analogous argument shows that
the same is true if the cycle of $\pi$ containing $k$ should have
an even length. Thus,
$$\det C=-{\sum _{}}{}^{\displaystyle\prime}
\underset k+\pi(k)\text { odd}\to{\sum _{k=2} ^{n}} (\sgn
\pi)\cdot a_{1,k}a_{1,\pi(k)} \underset i\ne k\to{\prod _{i=2}
^{n}}a_{i,\pi(i)},\tag\AL$$ where the sum is over all
$\pi\in\SS_{n-1}$ with only cycles of even length except that the
cycle containing $k$ has odd length.

To show that the expression (\AL) is equal to
$-2\Pf_e(A)\Pf_o(A)$, we construct a bijection between $\SS_{n-1}$
and $M_e\times M_o\times \{1,-1\}$, where $M_e$ denotes the set of
all perfect matchings on $\{1,2,\dots,n\}$ with the property that
$1$ is matched to an even number, and where $M_o$ is the analogous
set of perfect matchings with the property that  $1$ is matched to
an odd number. If $\pi$ is mapped to $(m_1,m_2,\et)$ under this
bijection, then this bijection will have the property that
$$\sgn \pi \cdot a_{1,k}a_{1,\pi(k)}
\underset i\ne k\to{\prod _{i=2} ^{n}}a_{i,\pi(i)}= -(\sgn
m_1)(\sgn m_2) \prod _{(i,j)\in m_1} ^{}a_{i,j}\prod _{(i,j)\in
m_2} ^{}a_{i,j}.$$ Clearly, given such a bijection, the assertion
of the theorem would be proved.

Let $\pi\in \SS_{n-1}$. Consider a cycle of $\pi$ not containing
$k$. Let $i$ be the smallest number in the cycle. Then we match
$i$ to $\pi(i)$ in $m_1$, $\pi(i)$ to $\pi^2(i)$ in $m_2$,
$\pi^2(i)$ to $\pi^3(i)$ in $m_1$, $\pi^3(i)$ to $\pi^4(i)$ in
$m_2$, etc. Considering the cycle containing $k$, we let $\et=1$
if $k$ is even while we let $\et=-1$ is $k$ is odd. If $k$ is
even, then we match $1$ to $k$ in $m_1$, $k$ to $\pi(k)$ in $m_2$,
$\pi(k)$ to $\pi^2(k)$ in $m_1$, etc., while if $k$ is odd, we
match $1$ to $k$ in $m_2$, $k$ to $\pi(k)$ in $m_1$, $\pi(k)$ to
$\pi^2(k)$ in $m_2$, etc. It is obvious that this sets up a
bijection. The fact that the sign behaves in the correct way under
the bijection can be shown in a similar manner as in the proof of
Proposition~2.2 in \cite{\StemAE}. \qed
\enddemo

\enddocument